\documentclass[12pt,a4paper]{article}
\usepackage[centertags]{amsmath}
\usepackage{amsfonts}
\usepackage{amssymb}
\usepackage{amsthm}
\usepackage{newlfont}
\usepackage[T1]{fontenc}
\usepackage[ansinew]{inputenc}
\usepackage{amsmath}
\usepackage{graphicx}
\usepackage{color}
\usepackage[ansinew]{inputenc}
\usepackage[all]{xy}
\theoremstyle{plain}
\newtheorem{thm}{Theorem}[section]
\newtheorem{cor}[thm]{Corollary}
\newtheorem{lem}[thm]{Lemma}

\theoremstyle{definition}
\newtheorem{defn}[thm]{Definition}

\newcommand{\al}{\alpha}
\newcommand{\be}{\beta}
\newcommand{\f}{\varphi}

\newcommand{\de}{\delta}
\newcommand{\s}{\sigma}

\newcommand{\e}{\varepsilon}
\newcommand{\dd}{\partial}
\newcommand{\g}{\gamma}
\newcommand{\G}{\Gamma}
\newcommand{\Si}{\Sigma}

\newcommand{\R}{\mathbb R}
\newcommand{\Z}{\mathbb Z}

\newcommand{\rel}{,}

\newcommand{\im}{\textup{Im}}
\newcommand{\Diff}{\textup{Diff}}
\newcommand{\Homeo}{\textup{Top}}
\newcommand{\iso}{\cong}
\newcommand{\he}{\simeq}

\newcommand{\set}[1]{\left\{#1\right\}}

\newcommand{\To}{\longrightarrow}
\newcommand{\into}{\hookrightarrow}

\newcommand{\abs}[1]{\left\vert#1\right\vert}

\newcommand{\id}{\textup{id}}
\newcommand{\pr}{\textrm{pr}}

\newcommand{\del}{\subseteq}
\newcommand{\fra}{\setminus}

\catcode`\=13
\def{$\bowtie$}
\author{Søren K. Boldsen}
\title{Different versions of \\mapping class groups of surfaces}
\date{\today}

\begin{document}
\maketitle

\section{Introduction}
Let $F$ be a compact connected smooth surface, possibly with boundary and not necessarily oriented. The objects of study in this paper are
\begin{eqnarray*}
\Diff(F,\set{\dd F})&=& \set{\f:(F,\dd F)\To (F,\dd F) \mid \f \text{ is a diffeomorphism}},\\
\Homeo(F,\set{\dd F})&=& \set{\f:(F,\dd F)\To (F,\dd F) \mid \f \text{ is a homeomorphism}},\\
\textup{hAut}(F,\set{\dd F})&=& \set{\f:(F,\dd F)\To (F,\dd F) \mid \f \text{ is a homotopy equivalence}}.
\end{eqnarray*}
The main theorem of this paper is the following:
\begin{thm}\label{t:main}Let $F$ be a compact surface and not
a sphere, a disk, a cylinder, a Möbius band, a torus, a Klein bottle, or $\R P^2$. Then
\begin{equation*}
\pi_0(\Diff(F,\set{\dd F})) \stackrel{\iso}{\To} \pi_0(\Homeo(F,\set{\dd F}))
\stackrel{\iso}{\To} \pi_0(\textup{hAut}(F,\set{\dd F}))
\end{equation*}
are bijections.
\end{thm}

This result is far from new, but this paper will present a thorough and self-contained proof of the following bijection
\begin{equation}\label{e:bijection}
\pi_0(\Diff(F,\set{\dd F})) \stackrel{\iso}{\To} \pi_0(\textup{hAut}(F,\set{\dd F})).
\end{equation}
To get the Main Theorem from this result, we will use the result of \cite{Epstein} Thm 6.4 without proof.

We consider slightly different versions of the groups, where we assume $F$ is oriented in the last two groups:
\begin{eqnarray*}
\Diff(F \rel \dd F)&=& \set{\f\in \Diff(F,\set{\dd F}) \mid \f|_{\dd F}=\id},\\
\Diff_+(F,\set{\dd F})&=& \set{\f\in \Diff(F,\set{\dd F}) \mid \f\text{ is orientation-preserving}},\\
\Diff_+(F \rel \dd F)&=&  \Diff(F \rel \dd F)\cap \Diff_+(F,\set{\dd F}),
\end{eqnarray*}
and similar for $\Homeo$ and $\textup{hAut}$. By orientation-preserving we mean that the orientation class $[F,\dd F]\in H_2(F,\dd F)$ is preserved by $\f_*$. From the Main Theorem we easily deduce
\begin{thm}\label{t:main2}Let $F$ be a compact surface and not
a sphere, a disk, a cylinder, a Möbius band, a torus, a Klein bottle, or $\R P^2$. Then there are bijections
\begin{itemize}
  \item[$(1)$] $\pi_0(\Diff(F,\dd F)) \stackrel{\iso}{\To} \pi_0(\Homeo(F,\dd F))
\stackrel{\iso}{\To} \pi_0(\textup{hAut}(F,\dd F)),$
  \item[$(2)$] $\pi_0(\Diff_+(F,\set{\dd F})) \stackrel{\iso}{\To} \pi_0(\Homeo_+(F,\set{\dd F}))
\stackrel{\iso}{\To} \pi_0(\textup{hAut}_+(F,\set{\dd F})),$
  \item[$(3)$] $\pi_0(\Diff_+(F,\dd F)) \stackrel{\iso}{\To} \pi_0(\Homeo_+(F,\dd F))
\stackrel{\iso}{\To} \pi_0(\textup{hAut}_+(F,\dd F)).$
\end{itemize}
\end{thm}
The standard definition of the mapping class group of a surface $F$ is $\G(F)=\pi_0(\Diff_+(F\rel\dd F))$. The last part of Theorem \ref{t:main2} shows that it does not matter whether one considers diffeomorphisms, homeomorphisms, or even homotopy equivalences, when working in the mapping class group.
\newline\newline
It is a pleasure to thank Jørgen Tornehave for many fruitful discussions and help during my work on this paper.

\tableofcontents\newpage
\section{Preliminaries}
\begin{defn}An isotopy $\psi$ of $F$ is a path in $\Homeo(F,\set{\dd F})$,
i.e. $\psi:F\times I\To F$ is continuous map such that
$\psi_t=\psi(-,t):F\To F$ is a homeomorphism for all $t\in I$, and
we say that $\psi_0$ and $\psi_1$ are isotopic.

\noindent An isotopy is smooth if we can exchange homeomorphism
with diffeomorphism in the above. We then say that $\psi_0$ and $\psi_1$ are
smoothly isotopic.
\end{defn}

\begin{lem}\label{l:S^1}Let $f:S^1\To S^1$ an orientation preserving
diffeomorphism. Then $f$ is smoothly isotopic to the identity via a
smooth isotopy $f_t:S^1\times I\To S^1$ such that the function
$F:S^1\times I\To S^1\times I$ given by $F(z,t)=(f_t(z),t)$ is a
diffeomorphism, and
\begin{equation*}
    f_t(z)=\left\{
             \begin{array}{ll}
                f(z), & \hbox{for $0\le t< \e$,}\\
z, & \hbox{for $1-\e<t\le 1$.}
             \end{array}
           \right.
\end{equation*}
\end{lem}
\begin{proof} Since $f$ is smooth it defines a smooth
function $\tilde{f}:\R\To \R$ by lifting $f$ under the universal
covering $\exp:\R\To S^1$. Now take a smooth bump function
$\rho:I\To I$ satisfying
\begin{equation*}
    \rho(t)=\left\{
              \begin{array}{ll}
                0, & t\le \e, \\
                1, & t\ge 1-\e.
              \end{array}
            \right.
\end{equation*}
Let $\tilde F:\R\times I\To \R$ be given by
$\tilde{F}(x,t)=\rho(t)\tilde{f}(x)+(1-\rho(t))x$. This now defines
an isotopy from $\tilde{f}$ to the identity, and
$F(\exp(x),t)=(\exp(\tilde{F}(x,t)),t)$ is a diffeomorphism.
\end{proof}

The idea of the following proof is due to J. Alexander.

\begin{lem}\label{l:disk}Let $D$ be a disk and $N$ a collar neighborhood of the boundary. Suppose $f:D\To D$ is a homotopy
equivalence which restricts to an orientation preserving
diffeomorphism of $N$ of the form $f(z,t)=(f(z),t)$ for $(z,t)\in
N$. Then $f$ is homotopic to a diffeomorphism relative to a smaller
collar neighborhood.
\end{lem}

\begin{proof}We can assume $f:D\To D$, where $D=\set{z\in\R^2\mid \abs{z}\le
1+\e}$, and $N=\set{z\in D\mid 1-\e<\abs{z}\le
1+\e}$. The tubular coordinates on $N$ are $s\in [0,2\pi]$
and $t\in (-\e, \e]$. We first construct a
homotopy $\f_x$, $x\in [0,1]$, which is constant in $x$ outside $N$, from $f$ to a
function $g$ such that $g(s,0)=(s,0)$ in tubular coordinates. We use
the isotopy $f_x(s)$ from Lemma \ref{l:S^1}, and set
\begin{equation*}
\textstyle \f_x(s,t)=(f_{x(1-\frac{1}{\e}|t|)}(s),t), \quad t\in
(-\e, \e]
\end{equation*}
in tubular coordinates. Then $\f_0=f$ and $\f_1(s,0)=(s,0)$, and
$\f_x$ is the identity on a collar neighborhood of $\dd D$ by Lemma
\ref{l:S^1}.

We now make a homotopy $\psi_x$, $x\in [0,1]$, from $g$ to the function $h$
satisfying
\begin{equation*}
    h(z)=\left\{
           \begin{array}{ll}
             g(z), & \hbox{$\abs{z}>1$;} \\
             z, & \hbox{$\abs{z}\le 1$.}
           \end{array}
         \right.
\end{equation*}
Let $D'=\set{z\in D\mid\abs{z}\le 1}$, and define the solid cone
\begin{equation*}
    C=\set{(z,x)\del D\times I|\abs{z}\le 1-x}
\end{equation*}
with bottom $D'\times \set{0}$ and top $(0,1)$, and set
\begin{equation*}
    \psi_x(z)=\left\{
             \begin{array}{ll}
               (1-x)f(\frac{z}{1-x}), & (z,x)\in C, \\
               z, & (z,x)\in (D'\times I)\fra C, \\
                g(z), & (z,x)\in (D\fra D')\times I.
             \end{array}
           \right.
\end{equation*}
This is clearly continuous and constitutes a homotopy from $g$ to
$h$ through maps which are the identity on a collar neighborhood
of $\dd D$, since $g$ is. We claim $h:D\To D$ is a diffeomorphism.
Clearly, $h: D'\To D'$ is a diffeomorphism, and by Lemma
\ref{l:S^1}, $h$ is smooth on $D$, and for $\abs{z}>1$, $h=g$ is a
diffeomorphism $D\fra D'\To D\fra D'$.
\end{proof}

A result we will use repeatedly is the following smooth version of the Schönflies curve theorem.

\begin{lem}\label{l:disk}Let $f:S^1\To F$ be a smoothly embedded simple closed curve homotopic to zero in a surface $F$. Then the closure of the interior of $f(S^1)$ is a smoothly embedded disk in $F$.
\end{lem}

\begin{proof} By Thm 1.7 in \cite{Epstein} we know that $f$ separates $F$ into two components, and that one of them (call it $D'$) is homeomorphic to a disk $D^2$. Thus $D'$ is a connected orientated smooth 2-manifold with 1 boundary component and with Euler characteristic $\chi(D')=1$. Now by the classification of smooth surfaces, $D'$ is a smooth disk.
\end{proof}

\begin{defn}
Let $\al$ be a smoothly embedded 1-submanifold in a surface $F$. By the surface cut up along $\al$, denoted $F\fra \al$, we will mean the surface with boundary $F\fra N(\al)$, where $N(\al)$ is a tubular neighborhood of $\al$ in $F$.
\end{defn}

\begin{lem}\label{l:nonsep} Let $\al: (I,\dd I)\To (F,\dd F)$ be a simple curve in a
surface $F$. If the cut-up surface $F\fra \al(I)$ is disconnected, then the induced map
$\al_*: H_1(I,\dd I)\To H_1(F,\dd F)$ is the zero map.
\end{lem}
\begin{proof}Let $\bar{\al}=\al(I)\del F$, and consider the long exact sequence for the triple $(\dd F, \bar\al \cup\dd F,F)$:
\begin{equation*}
 \xymatrix{H_1(\bar\al\cup\dd F,\dd F)\ar[r]^{\quad i_*}& H_1(F,\dd F)
\ar[r]^{j_*\quad}&  H_1(F,\bar\al\cup\dd F)\ar[r]&
H_0(\bar\al\cup\dd F,\dd F)}
\end{equation*}
Here $H_0(\al\cup\dd F,\dd F)=0$, so $j_*$ is surjective. Also
$H_1(F,\dd F)\iso \Z^{2g+r-1}$ for $F=F_{g,r}$. Since $F\fra \bar
\al$ is not connected, we can write $F\fra\bar\al=F_1\sqcup F_2$,
and by excision,
\begin{eqnarray*}
  H_1(F,\bar\al\cup\dd F) &\iso& H_1(F_1\sqcup F_2,\dd F_1\sqcup\dd F_2)\iso H_1(F_1,\dd F_1)\oplus H_1(F_2,\dd F_2)\\
   &\iso&\Z^{2g_1+r_1-1}\oplus \Z^{2g_2+r_2-1}.
\end{eqnarray*}
Here $g=g_1+g_2$ and $r+1=r_1+r_2$, so since $j_*$ is surjective, we
conclude that $j_*$ is an isomorphism. Thus $i_*=0$, and the following diagram shows that $\al_*=0$:
\begin{equation*}
 \xymatrix{H_1(I,\dd I)\ar[d]_{\al_*}\ar[r]^{\al_*}& H_1(F,\dd F)\\
H_1(\bar\al\cup\dd F,\dd F)\ar[ru]^{i_*}&}
\end{equation*}
\end{proof}

\section{Surjectivity}
In this section we will prove that the map in \eqref{e:bijection} is surjective, i.e. a homotopy equivalence of a surface $F$ is homotopic to a diffeomorphism. We first prove this for surfaces with non-empty boundary, and then use this to obtain the proof for closed surfaces. The result for surfaces with non-empty boundary is strongly inspired by \cite{Hempel}.

\begin{thm}\label{t:surj}Let $F$ and $G$ be compact surfaces with
non-empty boundaries. Suppose $\pi_1(F)$ is non-trivial. Let $f:(F,\dd
F)\To (G,\dd G)$ be a map such that $f_*:\pi_1(F)\To \pi_1(G)$ is
injective and $f|_{\dd F}:\dd F\To \dd G$ is a smooth embedding.
Then there is a homotopy $f_t:(F,\dd F)\To (G,\dd G)$ with $f_0=f$
and $ f_1:F\To G$ a diffeomorphism.
\end{thm}

\begin{proof}First consider each boundary component $J$ of $F$, and $K$ of
$G$ where $f(J)\del K$. We can assume each $J$ and $K$ has a collar neighborhood of the form $J\times [0,\e]$ and $K\times [0,\e]$, where the map $f$ has the form $f(x,t)=(f|_{J}(x),t)$, by gluing on small cylinders, extending $f$ as desired, and smoothing out. Since $f$ is continuous, it is homotopic to a smooth map, and we can choose the homotopy to be constant on the collar neighborhoods, so we can assume that $f$ is smooth  an embedding on a neighborhood of $\dd F$.

We are going to cut up $G$ by a non-separating arc $\al$ (i.e. an embedded
connected 1-manifold with boundary) connecting two boundary
components of $G$ in the image of $f$. We would like to cut up $F$
by $f^{-1}(\al)$. To do this we must ensure that $f^{-1}(\al)$ is
also an embedded 1-manifold. This holds if $f$ is transverse to
$\al$. By Thom's transversality theorem, $f$ can be
approximated by a smooth map $g$ transverse to $\al$ arbitrarily close to
$f$. Even better, $g$ can be chosen such that $g|_A=f|_A$ for a
closed subset $A\del F$ on which the transversality condition on $f$ is already satisfied. If we choose the arc $\al$ to have the form $\al=(x_0,t)$, $t\in
[0,\e]$ on the collars $K\times [0,\e]$ for some $x_0\in K$, then
clearly we can take $A=\bigcup_{J\in\pi_0(F)}J\times [0,\e]$ in the above. Since the transverse map $g$ is arbitrarily close to $f$, they are homotopic, and we can assume $f$ is transverse to $\al$.

Since $f|_{\dd F}:\dd F\To \dd G$ is an embedding we can see that
$f^{-1}(\alpha)$ must consist of one arc in $F$ and possibly a
number of embedded circles, and as $F$ is compact, there is a finite
number of circles. Since $f_*$ is injective, the circles must be
null-homotopic in $F$, thus they must each bound a disk $D_0$ in $F$.
Taking a slightly larger disk $D\supseteq D_0$, then $f(\dd D)$ must
be contained in a tubular neighborhood of $\alpha$. Since $\dd D$
is disjoint from $f^{-1}(\al)$, all of $f(\dd D)$ is to the same
side of $\al$ in the tubular neighborhood.

Now $D$ is a disk  and $f(\dd D)$ is contained in a disk $E\del G$ on one side of $\al$ in the tubular neighborhood. Thus we can make a map $h:D \To G$ with $h(D)\del E$ and such that $f|_{\dd D}= h|_{\dd D}$. This
gives a map $H: S^2\To G$ by mapping the lower hemisphere by $f$
and the upper hemisphere by $h$. Since $G$ is not $S^2$ or $\R P^2$, we
know $\pi_2(G)=0$, so the map $H$ can be extended to a map $D³\To G$, thus giving a homotopy from $f$ to $h$. This will reduce the number of circles in the preimage, and we can thus assume that $f^{-1}(\alpha)$ is just an arc in $F$. By transversality we can assume that we have a tubular neighborhood of $f^{-1}(\alpha)$ mapping to a tubular neighborhood of $\alpha$.

We can now cut $F$ along $f^{-1}(\alpha)$ and $G$ along $\alpha$, to
obtain $\hat{F}$ and $\hat{G}$. After cutting up $F$ and $G$ along
an arc, we will actually have manifolds with corners, $\hat{F}$ and
$\hat{G}$. But clearly we can smooth out these corners inside the
collar neighborhoods where $f:\hat{F}\To \hat{G}$ is smooth.

Now we would like to show that the process will not
separate $F$. Consider the situation when we cut up along a non-separating arc
$\al$ in $G$. We can parametrise $\al$ and think of it as a function
$\alpha: (I,\dd I) \To (G,\dd G)$. This induces a map $\alpha_*:
H_1(I,\dd I;\mathbb{Z}_2) \To H_1(G,\dd G;\mathbb{Z}_2)$. The
condition that $\alpha$ is nonseparating translates as $\alpha_*\neq
0$. By the above we can assume that $f^{-1}(\al)$ is a single arc, which we parametrize as $\tilde{\alpha}: (I,\dd I) \To (F,\dd F)$:
\begin{equation*}
\xymatrix{&(I,\dd I)\ar[dl]_{\tilde{\alpha}}\ar[dr]^{\alpha} & \\
            (F,\dd F) \ar[rr]^f &       &          (G,\dd G)}
\end{equation*}
On homology this induces the commutative diagram
\begin{equation*}
\xymatrix{&H_1(I,\dd I;\mathbb{Z}_2)\ar[dl]_{\tilde{\alpha}_*}\ar[dr]^{\alpha_*} & \\
            H_1(F,\dd F;\mathbb{Z}_2) \ar[rr]^{f_*} &       &          H_1(G,\dd
G;\mathbb{Z}_2)}
\end{equation*}
But since $\alpha_*\neq 0$ we get $\tilde{\alpha}_*\neq 0$ and
thus by Lemma \ref{l:nonsep}, $\tilde{\alpha}\del F$ is
nonseparating.

Now we show that $f_*:\pi_1(\hat{F})\to \pi_1(\hat{G})$
is still injective after cutting up. We use that $F$ is homotopic to $\hat{F}\cup I$, where $I$ is a small interval connecting two points
$b_0, b_1\in\dd \hat{F}$. Using that $\hat{F}$ is connected we
choose a path $J$ in $\hat{F}$ from $b_0$ to $b_1$, such that $I\cup
J$ form a loop. Now $F\he \hat{F}\vee S^1$ (by contracting $J$ in
$\hat{F}$ to a point). Then $i_*: \pi_1(\hat{F})\To \pi_1(F)$ is injective, since
$i_*:\pi_1(\hat{F})\To \pi_1(F)=\pi_1(\hat{F})*\mathbb Z$ is just
the inclusion in the first factor by van Kampen's theorem. Now it
follows from the commutative diagram
\begin{equation*}
\xymatrix{\pi_1(\hat{F})\ar@{^(->}[r]^{i_*}\ar[d]^{\hat{f}_*}       &    \pi_1(F)\ar@{^(->}[d]^{f_*}\\
                  \pi_1(\hat{G})\ar[r]^{i_*}           &    \pi_1(G)}
\end{equation*}
that $\hat{f}_* : \pi_1(\hat{F})\To \pi_1(\hat{G})$ is also
injective.

It remains to show that by cutting up $F$ and $G$ they have
to become disks at the same time. Firstly if $G$ is a disk, then
$f_*:\pi_1(F)\To \set{1}$ is injective, so $\pi_1(F)=\set{1}$, and this
implies that $F$ is also a disk (since $F$ is a surface with
boundary). Conversely, if $G$ is not a disk then
neither is $F$, since given a non-separating arc $\al$ in $G$ we have
shown above that there exists a non-separating arc in $F$.

We are down to the case where $f$ is a map from a disk to
a disk that is smooth in a collar of the boundary, and this case is
handled by Lemma \ref{l:disk}. We can glue the resulting smooth maps
on the pieces together again, since the collar neighborhoods of
the boundary of each piece (where the map is smooth) are fixed by the
homotopy from Lemma \ref{l:disk}. So we are done.
\end{proof}

\begin{cor}\label{c:surj}
Let $F$ and $G$ be compact surfaces with
non-empty boundaries. Suppose $\pi_1(F)$ is non-trivial. Let $f:(F,\dd
F)\To (G,\dd G)$ be a map such that $f_*:\pi_1(F)\To \pi_1(G)$ is
injective and $f|_{N(\dd F)}:N(\dd F)\To N(\dd G)$ is a smooth embedding, where $N(-)$ denotes a neighborhood. Then there is a homotopy $f_t:(F,\dd F)\To (G,\dd G)$ with $f_0=f$ and $ f_1:F\To G$ a diffeomorphism, such that $f_t=f_0$ on a neighborhood of $\dd F$.
\end{cor}

\begin{proof} Use the proof above, but skip the first part which proves that $f|_{N(\dd F)}:N(\dd F)\To N(\dd G)$ can be made into a smooth embedding.
\end{proof}

\begin{lem}\label{l:cylinder}Let $f_0,f_1:S^1\To F$ be disjoint
non-trivial two-sided embeddings in the surface $F$. Assume there exist $m,n\in \Z_{+}$
such that $f_0^n$ and $f_1^m$ represent the same free homotopy class
in $F$. Then there is an embedding $\f:S^1\times I\To F$ such that
$\f|_{S^1\times\set{i}}=f_i$ for $i=0,1$, so $f_0$ and $f_1$ bound
a cylinder.
\end{lem}
\begin{proof}This is a special case of \cite{Epstein}, Lemma 2.4.

We start by cutting $F$ up along $f_0$ and then gluing a disk onto
each of the two new boundary components; let $M$ be the connected component containing $f_1$ in the resulting surface. Since $f_0$ is null-homotopic in $M$, then so is $f_0^n$ and thus $f_1^m$. Now we will show that $f_1$ is null-homotopic in $M$, so that it bounds a disk in $M$. First if $\dd M\ne\emptyset$, then $\pi_1(M)$ is a free group and thus if $f_1^m=1$ then $f_1=1$. Else $\pi_1(M)$ is a free group modulo the relation $\dd =
\Pi_{i=1}^g[a_i,b_i]\in\pi_1(M)$ (oriented case) or
$\dd=\Pi_{i=1}^ga_i^2\in\pi_1(M)$ (unoriented case). If $f_1^m=1$ but $f_1\ne1$, $\pi_1(M)$ will have torsion, and by \cite{LySc}
Prop. 5.18, the only case that allows for torsion is the unoriented
case with $g=1$. Then the component of $M$
containing $f_1$ is an $\R P^2$, but then there are no non-trivial
two-sided embeddings of $S^1$. So there can be no torsion, and
$f_1=1$ in $\pi_1(M)$.

The disk in $M$ bounded by $f_1$ contains either one or two of the disks
glued onto $f_0$ to form $M$, since $f_1$ was non-trivial in $F$. If
the disk bounded by $f_1$ in $M$ contains just one glued-on disk,
then $f_0$ and $f_1$ together bound a disk blown up at one point; a
cylinder in $F$. In particular, if $f_0$ is separating, then the disk bounded by $f_1$ in $M$ contains just one glued-on disk, so we are done. Now if the disk bounded by $f_1$ in $M$ contains two of the
glued-on disks, then $f_1$ was separating in $F$, since we obtain
$F$ from $M$ by removing the glued-on disks and gluing up along
their boundaries. The cylinder can thus be obtained if we interchange
$f_0$ and $f_1$.
\end{proof}

The condition in the preceding Theorem \ref{t:surj} about the map $f$ being an
embedding on the boundary is not essential if $f$ is a homotopy
equivalence, as we show next:

\begin{lem}\label{l:notemb}
Suppose $f:(F,\dd F)\To (G,\dd G)$ induces an isomorphism
$f_*:\pi_1(F)\To \pi_1(G)$, and suppose $F$ is compact with $\dd F\neq \emptyset$ and is neither a disk, a cylinder nor a Möbius band. Then the following holds:
\begin{itemize}
  \item[$(i)$]For all boundary components $J\del  \dd F$ and $K\del
  \dd G$ such that $f(J)\del K$, the composite $\Z\iso \pi_1(J)\stackrel{f}{\To}
  \pi_1(K)\iso \Z$ is multiplication by $\pm1$, and no two
  different boundary components in $F$ are taken to the same boundary
  component in $G$.
  \item[$(ii)$]$f$ is homotopic to a map $g:(F,\dd F)\To (G,\dd G)$ with $g|_{\dd F}:\dd F\To \dd G$
an embedding.
\end{itemize}
\end{lem}

\begin{proof}
Let $J\del F$ be a boundary component, and let $K\del G$ be the
boundary component with $f(J)\del K$. We have a commutative diagram,
\begin{equation}\label{e:notemb}
\xymatrix{
\pi_1(J)\ar@{^(->}[r]_{(f|_J)_*}\ar@{^(->}[d] & \pi_1(K)\ar[d] \\
\pi_1(F)\ar[r]^{\iso}_{f_*}& \pi_1(G) }
\end{equation}
Here, the vertical map $\pi_1(J)\To \pi_1(F)$ is injective, since it is a non-zero map (as $F$ is not a disk) from $\pi_1(J)\iso\Z$ into the free group
$\pi_1(F)$. Then $(f|_J)_*$ is multiplication by an integer $n\ne 0$.

If $F$ has more than 1 boundary component, we can choose generators
for $\pi_1(F)$ such that the generator of $\pi_1(J)$ goes to a
generator of $\pi_1(F)$ under the left vertical map in
\eqref{e:notemb}. Since $f_*$ is an isomorphism, it takes generators
to generators, and thus it follows by commutativity that $n=\pm 1$.

If $F$ only has the one boundary component $J$, then the generator
$\alpha$ of $\pi_1(J)$ maps to either $\dd =
\Pi_{i=1}^g[a_i,b_i]\in\pi_1(F)$ (oriented case) or
$\dd=\Pi_{i=1}^ga_i^2\in\pi_1(F)$ (unoriented case). If
$f_*(\alpha)= x^n$ for a generator $x$ of $\pi_1(K)$, we get by
commutativity that $f_*(\dd)\in \pi_1(G)$ would be an $n$th power of
something. Since $f_*:\pi_1(F)\To \pi_1(G)$ is an isomorphism, $\dd$
itself would be an $n$th power of some element. In case $\dd=
a_1^2$, $F$ is a Möbius band, so this cannot happen. In all other
cases we get $n=\pm 1$.

We have shown that $(f|_J)_*:\pi_1(J)\To \pi_1(K)$ is an
isomorphism. Thus we can homotope $f$ in a collar neighborhood
around $J$ such that $f|_{J}:J\To K$ is a diffeomorphism. We do this
for every boundary component of $F$.

 All that is left is to check that no two boundary components $J_1, J_2$ of $F$ map to the same boundary component $K$ in $G$. If that were the case,
 the elements of $\pi_1(F)$ generating $\pi_1(J_1)$ and $\pi_1(J_2)$
 would both map to a generator of $\pi_1(K)$, i.e. would coincide up to
 a sign, since $f_*: \pi_1(F)\To \pi_1(G)$ is an isomorphism. Then by
Lemma \ref{l:cylinder}, $F$ would be a cylinder, which it is not.
\end{proof}

\begin{thm}\label{t:surjdone}Let $F$ and $G$ be compact surfaces, and let $f:F\To G$ be a homotopy equivalence. Assume neither $F$ nor $G$ is a disk, a sphere, a cylinder, a Möbius band, a torus, a Klein bottle, or $\R P^2$. Then $f$ is homotopic to a diffeomorphism.
\end{thm}

\begin{proof}If $F$ and $G$ have non-empty boundary, Lemma \ref{l:notemb} and Theorem \ref{t:surj} give the result. So assume that $F$ and $G$ are closed surfaces.

Let $B\del G$ be a non-separating, 2-sided simple closed curve in
$G$. Since $f$ is homotopic to a smooth map which is transverse to $B$, we can
assume that $f$ is smooth and transverse to $B$. Consider the components of $f^{-1}(B)$. By transversality and compactness, this is a finite set of disjoint 1-submanifolds of $F$. As in the proof of Theorem \ref{t:surj}, we can homotope $f$ so that no component in $f^{-1}(B)$ bounds a disk. For any 1-sided
simple closed curve $\g$ in $f^{-1}(B)$, take a small tubular
neighborhood $M$ of $\g$ such that $f(M)\del N$, where $N$ is a
tubular neighborhood of $B$. Since $M\fra \g$ is connected and
$f(M\fra \g)\del N\fra B$, it follows that $M\fra\g$ maps to the
same side of the 2-sided curve $B$ under $f$. This implies that we
can homotope $f$ in $M$ to a function not hitting $B$. So we can
assume that no component of $f^{-1}(B)$ is a 1-sided simple closed
curve.

Now let $H_0,H_1$ be two components of $f^{-1}(B)$, and let
$h_0,h_1:S^1\To F$ be parametrizations of $H_0$ and $H_1$,
respectively. Then
\begin{equation*}
\Z\iso \pi_1(H_i)\stackrel{f}{\To} \pi_1(B)\iso \Z
\end{equation*}
is multiplication by some $m_i\in \Z$. Note that $m_i\ne 0$ since
$h_i$ is nontrivial in $F$ and $f$ is injective on $\pi_1(F)$. This
gives that $f_*(h_0^{m_1})=f_*(h_1^{m_0})\in \pi_1(G)$, and since
$f$ is injective on $\pi_1(F)$, $h_0^{m_1}=h_1^{m_0}\in \pi_1(F)$.
Then by Lemma \ref{l:cylinder} they bound a cylinder (if $h_0$ and
$\bar{h_1}$ bound a cylinder then so do  $h_0$ and ${h_1}$). This
cylinder might contain components of $f^{-1}(B)$, but since there
are finitely many such components, we can take a cylinder whose
intersection with $f^{-1}(B)$ is precisely its ends, call them $h_0$
and $h_1$ again.

Now the cylinder gives a homotopy $c:S^1\times I\To
F$ from $h_0$ to $h_1$, and thus $f\circ c:S^1\times I\To G$ is a
homotopy in $G$, with $f(c(S^1\times ]0,1[))\cap B=\emptyset$. Thus
we get a continuous map $\widetilde{f\circ c}:S^1\times I\To G\fra
B$ into the cut-up surface $G\fra B$. This is a homotopy between
non-zero powers of boundary components of $G\fra B$. Now by Lemma
\ref{l:cylinder}, if these two boundary components are distinct,
$G\fra B$ would be a cylinder. But this is impossible, since $G$ is
neither a torus nor a Klein bottle. This implies that both ends of
the cylinder is mapped to the same boundary component in $G\fra B$,
and thus we can change $f$ by a homotopy to remove $h_0$ and $h_1$
from $f^{-1}(B)$ without changing $f^{-1}(B)$ otherwise. We can now
assume that $f^{-1}(B)$ is a single closed curve, since
$f^{-1}(B)=\emptyset$ implies that $f_*:\pi_1(F)\To \pi_1(G)$
factors through $\pi_1(G\fra B)$ but $\pi_1(G\fra B)\To\pi_1(G)$ is
not surjective. We can finally see that the curve $f^{-1}(B)$ is
non-separating by Lemma \ref{l:nonsep}, since $B$ is non-separating and $f_*:H_1(F)\To
H_1(G)$ is a group homomorphism.

Consider $f|:N(f^{-1}(B))\To N(B)$, where $N(-)$ denotes a
tubular neighborhood. Then, using a method as in the proof of Lemma \ref{l:S^1} on $f^{-1}(B)$ and a bump function to extend to $N(f^{-1}(B))$, one can see that $f$ is homotopic to a map $g$ with $g^{-1}(B)=f^{-1}(B)$, such that $g|_{N(g^{-1}(B))}$ is a smooth covering map (the number of sheets will be the degree of $f:f^{-1}(B)\To B$). So now we assume that $f$ is a smooth covering map on a neighborhood of $f^{-1}(B)$.

Since $f_*: \pi_1(F\fra f^{-1}(B))\To \pi_1(G\fra B)$ is
injective (\cite{LySc} prop 5.1), we can choose a covering $\rho:\widetilde{G\fra B}\To G\fra B$ and lift $f$ as in the diagram,
\begin{equation}\label{e:cover}
    \xymatrix{ & \widetilde{G\fra B}\ar[d]^{\rho} \\
            F\fra f^{-1}(B)\ar[ur]^{\tilde f}\ar[r]^{\quad f}& G\fra B}
\end{equation}
such that $\tilde{f}_*:\pi_1(F\fra f^{-1}(B))\iso
\pi_1(\widetilde{G\fra B})$. Moreover, $\rho$ is a finite-sheet covering, since $f$ maps (a parametrization of) $f^{-1}(B)$ to a non-zero multiple of (a
parametrization of) $B$, and the number of sheets is locally
constant. So $\widetilde{G\fra B}$ is compact.

Now in a neighborhood of the boundary of $F\fra f^{-1}(B)$, $\tilde f$ is a covering map, and $\tilde f_*$ is an isomorphism on $\pi_1$. So $\tilde f$ is an embedding on a neighborhood of the boundary. By Corollary \ref{c:surj} on $\tilde f:F\fra f^{-1}(B)\To
\widetilde{G\fra B}$, $\tilde f$ is homotopic to a diffeomorphism,
relative to a neighborhood of the boundary. Glue up this diffeomorphism to a map $g:F\To G$ which will be homotopic to $f$, and be both a homotopy
equivalence and a smooth covering map. The last two imply that $g$ is a diffeomorphism $F\To G$.
\end{proof}

\section{Injectivity}
In this section we will prove that the map in \eqref{e:bijection} is injective, i.e. if a diffeomorphism is homotopic to the identity, it is smoothly isotopic to the identity.

\begin{defn}Let $f,g:I\To F$ be smooth embeddings into a surface $F$. We say
that $f$ and $g$ \emph{form an ''eye''} if the following is
satisfied:
\begin{itemize}
  \item[$(i)$]$f(I)\cup g(I)$ bounds a disk in $F$.
  \item[$(ii)$]$f|_{[0,\e[} = g|_{[0,\e[}$, $f|_{]1-\e,1]} = g|_{]1-\e,1]}$, and
$f$ is disjoint from $g$ on $]\e,1-\e[$.
\end{itemize}
\end{defn}
\begin{lem}\label{l:eye}
Let $f,g: I\To F$ be two smooth embeddings into a surface $F$ which
form an ''eye''. Then there is a smooth isotopy $\f_t$ of $F$ with $\f_0=\id_F$ and
$\f_1\circ g=f$. Furthermore, there is a small neighborhood $N$ of
the disk bounded by $f$ and $g$ for which $\f_t$ is the identity
outside $N$ for all $t$.
\end{lem}

\begin{proof}
Let $N_f$ be a tubular neighborhood of $f(I)$, given by a normal
vector field $\xi_f$ to $N_{f}$. Let also $N_g$ be a tubular
neighborhood of $g(I)$ given by a normal vector field $\xi_g$, such
that $N_f\cup N_g$ is an annulus. This is possible since $f(I)\cup
g(I)$ bounds a disk in $F$.

There is a diffeomorphism $\psi_f: N_f\To V_f\del\R^2$ such that
$\psi_f\circ f$ is the standard embedding $I\To \R\times \set{0}$.
We can take $V_f=I\times ]-\e,\e[$. We want to extend $\psi_f$ to a
diffeomorphism $\psi_{fg}:N_f\cup N_g\into \R^2$, i.e.
$\psi_{fg}|_{N_f}=\psi_f$.

First we note that inside $V_f=I\times ]-\e,\e[$ we have the image
\begin{equation*}
G = \psi_f(g(I)\cap N_f).
\end{equation*}
By taking $\e$ small, we can ensure that $G$ is the graph
$\set{(t,h(t))}$ of a smooth function $h:[0,\de[\cup]1-\de,1]\To
[0,\infty[$. We can extend $\psi_f$ to a map $\tilde{\psi}_{fg}$
defined on $N_f\cup g(I)$ such that $\tilde{\psi}_{fg}\circ g : I\To
\R^2$ is smooth, using bump functions etc as usual, such that the
image $G_{I}=\tilde{\psi}_{fg}\circ g(I)$ is the graph
$\set{(t,h(t)}$ of a function $h:I\To [0,\infty[$, see Figure 1.

\begin{figure}[htb!]\label{f:tub}
\begin{center}
\includegraphics[width=0.7\textwidth]{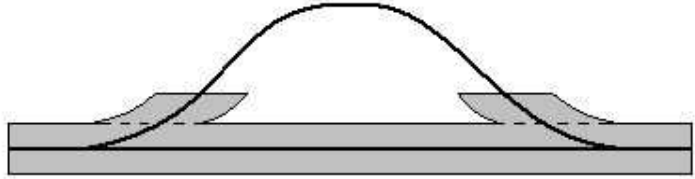}
\end{center}\caption{The tubular neighborhood $V_f$ and the graph $G_I$ of $h$ in $\R^2$.}
\end{figure}

We define a tubular neighborhood of $G$ using the vector field
$\eta_G=(\psi_{f})_*(\xi_g|_{N_f\cap N_g})$. Since $\psi_{f}$ is a
diffeomorphism, $\eta_G$ is a transverse vector field, and so
defines a tubular neighborhood $N_G$ of $G$ inside $V_f$. Now we
shrink $V_f$ to $I\times ]-\e',\e'[$ where $\e'<\e$ (thus also
shrinking $N_f$). Then we cover $G_I$ by two open sets in $\R^2$,
$U_1$ covering $G_{I}\fra G$, and $U_2$ whose intersection with $U_1$ lie in $N_G$ and outside $V_f$, see Figure 2. Then we take a partition of unity $\rho_1$, $\rho_2$ with respect to $U_1$,
$U_2$.
\begin{figure}[htb!]\label{f:tub2}
\begin{center}
\includegraphics[width=0.7\textwidth]{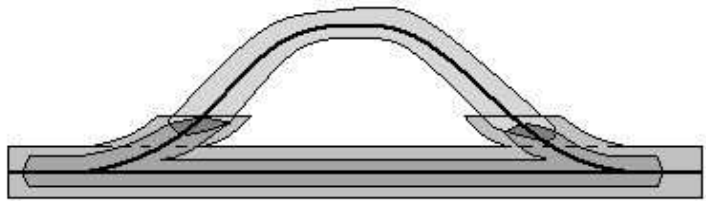}
\end{center}\caption{Neighborhoods $U_1$ and $U_2$ of $G_I$.}
\end{figure}

Let $\eta_{I}$ be the standard normal vector field to $G_{I}$,
defined on $G_{I}\fra G$. Then we make a new vector field
$\rho_1\eta_{I}+\rho_2\eta_G$. Since $\rho_1\eta_{I}+\rho_2\eta_G$
is never 0 or tangent to $G$, this defines a tubular neighborhood
$V_g$ of $G_{I}$. This tubular neighborhood coincides with $N_G$ on
$V_f$, and thus gives a diffeomorphism $\psi_{fg}:N_f\cup N_g \To V_f\cup V_g$ which extends $\tilde{\psi}_{fg}$.

The inner boundary circle $C$ of the annulus $N_f\cup N_g$ bounds a
disk $D'\del F$, and so the image $\psi_{fg}(C)$ also bounds a disk $D_{\R^2}\del \R^2$. Then we can extend $\psi_{fg}|_C$ to a map $D'\To D_{\R^2}$, which is necessarily a homotopy equivalence, so by Lemma
\ref{l:disk} we can replace it by a diffeomorphism $\psi_{D'}:D'\To
D_{\R^2}$, such that $\psi_{D'}|_C=\psi_{fg}|_C$. The we can glue
$\psi_{D'}$ and $\psi_{fg}$ along $C$ to obtain a diffeomorphism
$\Psi$ from $D=D'\cup N_f\cup N_g$ onto a disk in $\R^2$.

Now we can use a vertical flow in $D'\cup N_g\cup N_f$ (i.e. a pullback under $\Psi$ of the obvious vertical flow in $\R^2$) to make $\im(g)=\im(f)$, and lastly a horizontal flow in $N_f$ to make $g=f$.
\end{proof}

\begin{lem}\label{l:trans}Given two smoothly embedded arcs $f,g: I\To F$ satisfying $f(\set{0,1})
\cap g(I)=f(I) \cap g(\set{0,1})=\emptyset$. Then there is a smooth
isotopy $\f_t$ of $\id|_F$ such that $\f_1\circ f$ and $g$ intersect
transversally. Moreover $\f_t$ is the identity outside a tubular
neighborhood of $f$.
\end{lem}

\begin{proof}
Take an open tubular neighborhood of $f$, $N_f$, of constant radius,
where $r:N_f\To f(I)$ is the retraction. Inside ${N}_f$ take a
closed tubular neighborhood of $f$ of constant radius, $N^c_f$.
We cover $g(I)\cap N^c_f$ with sets of the form ${N}_f(a,b)=
\set{x\in {N}_f\mid f^{-1}(r(x))\in ]a,b[}$, where $a<b\in I$, and
$f(a), f(b)$ is outside $g(I)$. Since $g(I)\cap N^c_f$ is compact,
we can assume that it is a finite covering, ${N}_f(a_i,b_i)$,
$i=1,\ldots, N$, where $a_1<a_2<\cdots<a_N$. For each $x\in F$ where $f$ and $g$ intersect
non-transversally, $x\in{N}_f(a_i,b_i)$ for some $i$. Now take the
first such $i$. Then we can choose another arc $\tilde{g}: I\To
{N}_f(a_i,b_i)$ such that $g$ and $\tilde g$ form an ''eye'' and
$\tilde{g}$ and $f$ intersect transversally  for all $x\in \tilde
g(I)\cap f(I)\del {N}_f(a_i,b_i)$. Now by Lemma \ref{l:eye} there is
an isotopy from $g$ to $\tilde{g}$ in ${N}_f(a_i,b_i)$, which is the
identity outside ${N}_f(a_i,b_i)$. Doing this for each $i$, we
obtain in finitely many steps an isotopy which is the identity
outside $N_f$, making $f$ and $g$ intersect transversally.
\end{proof}

\begin{lem}\label{l:Baer}Let $F$ be a compact surface with $F\neq \R P^2$, $S^2$, and let $f$ be a
diffeomorphism of $F$.
\begin{itemize}
  \item[$(i)$]Let $\al_i:S^1\To F\fra \dd F$ be a finite family of disjoint,
  non-trivial, pairwise
  non-homotopic two-sided simple closed curves, with $f\circ \al_i \he \al_i$ for all $i$. Then
  there is an smooth isotopy $f_t$ of $F$ such that $f_0=f$ and $f_1\circ
  \al_i = \al_i$ and the identity extends to tubular neighborhoods.
  \item[$(ii)$]Let $\al_i:I\To F$ be a finite family of simple
  curves, disjoint except possibly at endpoints, with $f\circ
\al_i\he \al_i$ and $f\circ \al_i=\al_i$ near the endpoints for all $i$. Let $A\del F$ be a union of disjoint
non-trivial closed curves, with $f|_A=\id$ and $\al_i(I)\cap
A=\al_i(\dd I)$ for all $i$. Then there is an smooth isotopy $f_t$
of $F$, such that $f_0=f$, $f_1\circ \al_i=\al_i$ and the identity
extends to tubular neighborhoods. Furthermore $f_t|_A=\id$ for all
$t$.
\end{itemize}
\end{lem}
\begin{proof}
$(i)$ and $(ii)$ can be proved by the same methods, so we handle the
two cases as one initially. But we will also use $(i)$ to prove
$(ii)$. First, in both cases we have a closed subset $A\del F$ with
$f|_A=\id$ (in case $(i)$, $A$ starts as $\emptyset$). Consider a
single curve $\al=\al_1$. We will make an isotopy $f_t$ of $F$ such
that $f_0=f$, $f_1\circ \al=\al$, and $f_t|_A=\id$ for all $t$. Then
we can let $A_1= A\cup \al(I)$, and use the result for $f_1$ and
$A_1$ on $\al_2$, completing the proof in a finite number of steps. So consider a curve $\al$ as in $(i)$ or $(ii)$, and let $\be =
f\circ \al$ be the image curve. By assumption, $\be\he\al$.

In case $(ii)$, there are small neighborhoods $N_0$ and $N_1$ of the start and end points where $\al$ and $\be$ agree. Inside $N_0$ and $N_1$ we can make an isotopy of $f$ which perturbs $\be$ slightly, so that $\al$ and $\be$ agree near the start/end point, and then become disjoint. By shrinking $N_0$ and $N_1$ we can assume that $\al$ and $\be$ are disjoint on $\dd N_0$ and $\dd N_1$. Our goal is now to make $\al$ and $\be$ disjoint outside $N_0$ and $N_1$. From now on, we will ignore $N_0$ and $N_1$ in the proof, and only work with $\al$ and $\be$ outside them.

By Lemma \ref{l:trans} we can assume $\al$ and $\be$ are transverse to each other. Then $\al$ and $\be$ have finitely many intersection points by
compactness. To get an isotopy of $F$ taking $\be$ to $\al$, we
will first ensure that $\al$ and $\be$ have no intersection points. To do
this, consider the universal covering $\pi:\tilde F\To F$. We can
model $\tilde F$ as an open disk in $\R^2$. Take a fixed lift
$\tilde\be$ of $\be$.

We consider all the connected components of $\pi^{-1}(\al)$ that
intersect $\tilde \be$. There are finitely many such components,
call them $\tilde\al_k$, since $\al$ and $\be$ have finitely many
intersection points. The $\tilde\al_k$ are also transverse to
$\tilde\be$. Now we look for a pair of intersection points between
$\tilde \be$ and an $\tilde \al_i$, such that the part of the two
curves between these points (a closed curve, call it $\s$) bounds a
disk whose interior does not contain any points on $\tilde\be$ or
$\tilde\al_k$ for any $k$. So $\s$ is a simple closed curve in
$\tilde F$ bounding a disk. Projecting onto $F$, we get $\pi\circ\s$ (the parts of $\al$ and of $\be$ between two intersection points)
also a simple closed curve, which is null-homotopic, so according to
Lemma \ref{l:disk}, $\pi\circ\s$ bounds a disk in $F$. We can choose a curve $\be'$ which form an ''eye'' with $\be$ and which does not intersect $\al$ in a neighborhood of the disk bounded by $\pi\circ \s$. Then by lemma \ref{l:eye} we can isotope $\be$ to $\be'$, so that there are two fewer intersection points between $\al$ and $\be'$. Since there are finitely many intersection points, this procedure terminates.

But we must show why we can always find such a $\s$ in $\tilde F$.
Since $\tilde\al_k$ is a connected component of $\pi^{-1}(\al)$,
each $\tilde\al_k$ separates $\tilde F$. So if $\tilde \be$ crosses
$\tilde \al_i$ once, it must cross it again (let us choose the first
time it does so), as it is transverse to $\tilde\al_i$. Now $\tilde
F\del \R^2$, so the part of $\tilde\be$ and   $\tilde\al_i$ between
these two intersection points will bound a disk. If this disk
contains parts of $\tilde \be$ or $\tilde \al_k$'s, there will be a
smaller disk inside which satisfies the requirements, since there
are finitely many intersection points. In this way we can isotope
$\be$ to a curve which does not intersect $\al$ (in case $(ii)$,
except in $N_0$ and $N_1$).

In case $(i)$, we now have two homotopic disjoint simple closed
curves $\al$ and $\be$. Then according to  Lemma 1.4, they bound a
cylinder. Recall that the set $A$ (fixed by $f$) consists of the
curves already handled, i.e. a union of non-trivial closed curves,
none of which are homotopic to $\al$, and thus not to $\be$, either.
Thus $A$ cannot intersect the cylinder bounded by $\al$ and $\be$
(in fact, $A$ cannot intersect a small open neighborhood of the
cylinder). Then clearly there is an isotopy $f_t$ of $F$, which is
the identity on $A$, taking $\be$ to $\al$.

In case $(ii)$, the two curves $\al$ and $\be$ are homotopic and
form a simple closed curve, so again they bound a disk. Recall that $A$ originally consisted of non-trivial closed curves, so none of these can be
inside the disk. As we add curves to $A$, the circles get connected
by arcs. None of these can intersect $\al$, since they were assumed
to be disjoint from the start. As $f$ is the identity on $A$, they
cannot intersect $\be=f\circ \al$, either. Thus $A$ cannot cross the
boundary of the disk, so $A$ and the disk are disjoint. Thus by
lemma \ref{l:eye} we can make an isotopy $f_t$ of $F$, which is the
identity on $A$, so that $f_1\circ\al=\al$.

Now we extend the result to tubular neighborhoods of the curves. We
make a tubular neighborhood $M_0$ of $\al$, and by compactness
identify it with $S^1\times ]\!\!-\!\e,\e[$ in case ($i$) and $I\times
]\!\!-\!\e,\e[$ in case ($ii$). Now for $(x,t)$ in a smaller neighborhood
$M_1\subset M_0$ of $\al$, the projection the second coordinate $\pr_tf_x(t):=\pr_t(f(x,t))$ has positive differential, and thus for all $x$ the image of $f_x(t)$, $\set{(x',t')\mid (x',t')=f_x(t) \text{ for some } t\in ]\!\!-\!\e,\e[}$ is the graph of a function $h_x(t')=x'$. Now we can make tubular neighborhood $M_2$ such that $M_2\subset f(M_1)$ and by possibly shrinking it
assume that $M_2=I\times ]\!-\!\de,\de[$ or $M_2=S^1\times ]\!-\!\de,\de[$. For definiteness, say $M_2=I\times ]\!\!-\!\de,\de[$. Choose a smooth
bump function $\rho(t)$ with $\rho(t)=1$ for $|t|\le ½\de$ and
$\rho(t)=0$ for $|t|=\de$. Let
\begin{equation*}
  g_s(x,t)=\left\{\begin{array}{ll}
      \big((1\!\!-\!\!s)h_x(t')\!+\!s\big(\rho(t')x\!+\!(1\!\!-\!\!\rho(t'))h_x(t')\big), t'\big) & \hbox{for $(x,t)\!\in f^{-1}(M_2)$} \\
                                                                                f(x,t) & \hbox{otherwise.}
                                                                              \end{array}
                                                                            \right.
\end{equation*}
where $t'$ is the second coordinate of $f(x,t)$ as above. Then $g_s$
defines an isotopy from $f$ to a function $g_1$ with the property
that $g_1(x,t)=(x,t')$ for $\textstyle
t'\in]\!\!-\!\frac{\de}{2},\frac{\de}{2}[$. Now by stretching the
parameter $t'$ in each interval $\set{x}\times ]\!\!-\!\de,\de[$, we can
assume that $f$ is the identity on a (smaller)
neighborhood.
\end{proof}


\begin{cor}\label{c:Baer}If we in addition to the requirements in lemma \ref{l:Baer} require that $f$ is the identity on $\dd F$, then the isotopy can be assumed also to be the identity on $\dd F$.
\end{cor}
\begin{proof}All the steps in the proof can be done away from the boundary.
\end{proof}

\begin{thm}\label{t:inj}Let $F\neq S²,\R P²$ and let $f,g\in \Diff(F,\dd F)$ be homotopic. Then $f$ and $g$ are smoothly isotopic.
\end{thm}

To prove this I use the following result from \cite{Smale} without proof.

\begin{thm}[Smale]\label{t:smale} Let $f\in \Diff(D²,\dd D²)$. Then
  $f$ is smoothly isotopic to the identity, and if $f$ is the identity on the boundary then so is the isotopy.
\end{thm}
\begin{proof}[Proof of Theorem \ref{t:inj}] If we prove that $f^{-1}g$ is
  smoothly isotopic to the identity, we will have a smooth isotopy from $g$ to
  $f$. Thus we can restrict our attention to the case $f\he \id$.

  Choose a pair of pants/annular decomposition of the surface $F$, i.e. a
  collection of disjoint simple closed curves $\al_i:I\To F,
  i=1,\ldots,n$, in $F$. By Lemma \ref{l:Baer} $(i)$, $f$ is
  smoothly isotopic to a map $g$, which is the identity on a tubular neighborhood of the $\al_i$. In
  each pair of pants $P$, chose two curves that cut $P$ up into a
  disk (for each annulus, choose one curve). By Lemma \ref{l:Baer} $(ii)$, there is an isotopy of $F$,
  which is the identity on the $\al_i$, from $g$ to a map $h$ fixing
  a tubular neighborhood of the two curves in each pair of pants.
  Then we can use Smale's Theorem \ref{t:smale} on each disk,
  getting an isotopy to the identity.
  \end{proof}
\begin{cor}\label{c:inj}If we in addition to the requirements in theorem \ref{t:inj} require that $f$ and $g$ is the identity on $\dd F$, then the isotopy can be assumed also to be the identity on $\dd F$.
\end{cor}
\begin{proof} This is done as in theorem \ref{t:inj}, except that we use corollary \ref{c:Baer} instead of lemma \ref{l:Baer}, and using the fact that the isotopy in theorem \ref{t:smale} can be chosen to be the identity on boundary.
\end{proof}
\section{Proof of the Main Theorem}
As explained in the introduction, we will use a result of Epstein to prove the statement about $\Homeo(F,\set{\dd F})$:
\begin{thm}[Epstein]\label{t:injEpstein}Let $F$ a compact surface and let $f:F\To F$
be a homeomorphism homotopic to the identity. Then $f$ is isotopic to the identity.
\end{thm}
\begin{proof}This is a part of \cite{Epstein} Thm 6.4, which states exactly this result, but for maps preserving a basepoint. And clearly, by an
isotopy we can assume that $f$ preserves any given point $x_0$,
and then $f$ will be homotopic to the
identity through maps preserving $x_0$.
\end{proof}

Now we are ready to prove the bijections of the Main Theorem 1.1:
\begin{equation*}
\pi_0(\Diff(F,\set{\dd F})) \stackrel{\iso}{\To} \pi_0(\Homeo(F,\set{\dd F}))
\stackrel{\iso}{\To} \pi_0(\textup{hAut}(F,\set{\dd F}))
\end{equation*}

\begin{proof}[Proof of Theorem \ref{t:main}]Suppose $F$ is not
a sphere, a disk, a cylinder, a Möbius band, a torus, a Klein bottle, or $\R P^2$. Consider the composite map from \eqref{e:bijection},
\begin{equation}
\pi_0(\Diff(F,\set{\dd F})) \To \pi_0(\textup{hAut}(F,\set{\dd F})).
\end{equation}
According to Theorem \ref{t:surjdone}, the map is surjective, and by
Theorem \ref{t:inj}, it is injective. Now all that is left is to
show that
\begin{equation*}
\pi_0(\Homeo(F,\set{\dd F})) \To \pi_0(\textup{hAut}(F,\set{\dd F}))
\end{equation*}
is injective. But that is Theorem \ref{t:injEpstein}.
\end{proof}

We now deduce Theorem \ref{t:main2}:
\begin{proof}[Proof of Theorem \ref{t:main2}]Suppose $F$ is not
a sphere, a disk, a cylinder, a Möbius band, a torus, a Klein bottle, or $\R P^2$.

Similar to the proof of Theorem \ref{t:main}, we consider the composite
\begin{equation*}
\pi_0(\Diff(F\rel\dd F)) \To \pi_0(\textup{hAut}(F\rel\dd F)).
\end{equation*}
We can assume $\dd F\neq \emptyset$, otherwise this is the Main Theorem. By Cor. \ref{c:surj}, it is surjective, and by Cor. \ref{c:inj} it is injective.
To prove the result, it suffices to show that
\begin{equation}\label{e:mangler}
\pi_0(\Diff(F\rel\dd F)) \To \pi_0(\Homeo(F\rel\dd F))
\end{equation}
is surjective.
Consider the following fibration,
\begin{equation}\label{e:fibration}
\Diff(F\rel\dd F)) \To \Diff(F,\set{\dd F})\To \Diff(\dd F).
\end{equation}
Here, $\Diff(\dd F)$ is a semidirect product $\Si_n\ltimes \Diff(S^1)$, where $\Si_n$ denotes the symmetric group of permutations of $n$ elements, and $n$ is the number of boundary components of $F$. We have of course a similar fibration for $\Homeo$. We use $\Diff(S^1)\stackrel{\iso}{\To} \Homeo(S^1)$, and this implies
\begin{equation}\label{e:rand}
    \Diff(\dd F)\stackrel{\iso}{\To} \Homeo(\dd F).
\end{equation}
Now apply the long exact sequence of homotopy groups for the fibration \eqref{e:fibration} and its counterpart for $\Homeo$. Using \eqref{e:rand} and the Main Theorem, we get by the 5-lemma that the map \eqref{e:mangler} is surjective.

Now assume $F$ is oriented. We can write $\Diff(F,\set{\dd F})$ as the disjoint union
\begin{equation*}
    \Diff(F,\set{\dd F})=\Diff_+(F,\set{\dd F})\sqcup \Diff_-(F,\set{\dd F}),
\end{equation*}
where the latter denotes the orientation-reversing maps. Similarly for $\Homeo$ and $\textup{hAut}$. Since the maps in the Main Theorem respect this disjoint union, we immediately get the second part of \ref{t:main2}.

By the same argument we can deduce the last part of \ref{t:main2} from the first part.
\end{proof}


\begin{thebibliography}{20}
\bibitem[Epstein]{Epstein} D. B. A. Epstein  \emph{Curves on 2-manifolds and
isotopies}, Acta Mathematica 115 (1966), 83-107.
\bibitem[Hatcher]{Hatcher}A. Hatcher, \emph{Algebraic Topology}, Cambridge University press, 2002.
\bibitem[Hempel]{Hempel}J. Hempel, \emph{3-manifolds}, Ann. of Math. Princeton University press, 2002.
\bibitem[Lyndon-Schupp]{LySc}R. Lyndon and P. Schupp, \emph{Combinatorial Group Theory}, Springer, 1977.
\bibitem[Smale]{Smale}S. Smale, Diffeomorphisms of the 2-sphere, Proc. Amer. Math. Soc., 10 (1959), 621-626.
\end{thebibliography}
\end{document}